\documentclass[11pt,titlepage]{amsart}
\usepackage[]{atmcs}
\usepackage{amssymb,amsmath,latexsym,times,enumitem,tikz,soul, amsthm, float, graphicx,gensymb}

\usepackage[colorlinks=true, linkcolor=blue, citecolor=blue, urlcolor=blue]{hyperref}

\title{Topology of The Polar Vortex and Montana Weather}


\author[JD]{Joshua Dorrington}
\address{University of Bergen\\
    Bergen, Norway}
\email{joshua.dorrington@uib.no}

\author{Sushovan Majhi}
\address{George Washington University\\    
    Washington, D.C., USA}
\email{s.majhi@gwu.edu}

\author{Atish Mitra}
\address{Montanta Technological University\\    
    Montana, USA}
\email{AMitra@mtech.edu}

\author{James Moukheiber}
\address{The George Washington University, Washington D.C., USA.}
\email{jamesmoukh@gwu.edu}

\author{Demi Qin}
\address{Tulane University, New Orleans, USA.}
\email{yqin2@tulan.edu}

\author{Jacob Sriraman}
\address{Montanta Technological University\\    
    Montana, USA}
\email{jsriraman@mtech.edu}

\author[KS]{Kristian Strommen}
\address{University of Oxford\\
    Oxford, UK}
\email{kristian.strommen@physics.ox.ac.uk}

\begin{document}

\begin{abstract}
This paper explores the use of Topological Data Analysis (TDA) to investigate patterns in zonal-mean zonal winds of the Arctic, which make up the polar vortex, in order to better explain polar vortex dynamics. We demonstrate how TDA reveals significant topological features in this polar vortex data, and how they may relate these features to the collapse of the stratospheric vortex during the winter in the northern hemisphere. Using a time series representation of this data, we build a point cloud using the principles of Takens' Embedding theorem and apply persistent homology to uncover nontrivial topological structures that provide insight into the dynamical system's chaotic and periodic behaviors. These structures can offer new perspectives on the dynamics of the polar vortex, and perhaps other weather regimes, all of which have a global impact. Our results show clear transitions between seasons, with substantial increases in topological activity during periods of extreme cold. This is particularly evident in the historically strong polar vortex event of early 2016. Our analysis captures the persistence of topological features during such events and may even offer insights into vortex splitting, as indicated by the number of distinct persistent features.  This work highlights the potential of TDA in climate science, offering a novel approach to studying complex dynamical systems.
\end{abstract}

\maketitle

\section{Introduction}
Topological Data Analysis (TDA) has recently emerged as a powerful tool in data analysis, identifying otherwise invisible patterns in a noise-resistant manner \cite{carlsson_topology_2009}. TDA has been previously applied to explore climate and atmospheric science, such as atmospheric river patterns in large climate datasets, and the North Atlantic Jet Stream \cite{strommen_topological_2023, muszynski_topological_2019}.
In this paper, we propose using TDA to study the topological significance of the polar vortex of the northern hemisphere and draw conclusions on weather patterns in influenced regions.   



\subsection{Dynamical Systems and weather regimes}
We can categorize the local non-Gaussian nature of weather patterns, such as those caused by the shifting polar vortex, through \textit{weather regimes}. 
In essence, the goal of defining regimes is to identify a limited set of significant large-scale flow patterns that govern low-frequency variability and transition between one another in an approximately Markovian manner \cite{strommen_topological_2023}.
These weather regimes have been applied to multiple atmospheric patterns such as the North Atlantic Oscillation (NAO) \cite{hurrell_north_2003}. Simpler \textit{dynamical systems} like the Lorenz `63 (L63) system \cite{lorenz_deterministic_1963} were constructed that reflect similar characteristics while being much easier to simulate. 
A dynamical system is a system whose state is uniquely specified by a set of variables and whose behavior is described by predefined rules.  
In the case of weather regimes, we work with \textit{continuous} dynamical systems. 
Such systems may behave chaotically for some parameters, meaning that small changes in initial conditions can lead to deterministic yet drastically different outcomes \cite{lorenz_deterministic_1963}.

\subsection{Polar Vortices}
The \textit{polar vortex} is a large-scale cyclonic circulation of cold air above the polar regions of Earth, driven by temperature gradients between the Arctic (or Antarctic) and mid-latitudes. In the Northern Hemisphere, the polar vortex is strongest in winter as the gradient between polar and warmer mid-latitude air is at its greatest. \cite{manney_whats_2022} This vortex is crucial in shaping atmospheric dynamics, influencing jet streams, and globally impacting weather patterns. Sudden stratospheric warming events contribute to the instability of the polar vortex, which can trigger extreme cold weather in affected areas and shifts in mid-latitude climate. \cite{schoeberl_structure_1992}
In this paper, we see how analyzing the topological structure of the polar vortex and its wind patterns offer new insights into the formation, breakdown, and broader climate effects, and show a nontrivial link between the transitions of the polar vortex and the significant topological features obtained through time-delay embeddings. 

\begin{figure}[h]
    \centering
    \includegraphics[width=0.6\linewidth]{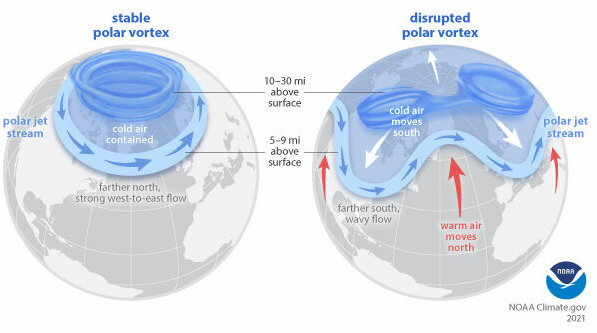}
    \caption{Left: A stable polar vortex. Right: Warm air from the south disrupts the vortex \cite{butler_understanding_2021}.}
    \label{fig:polar-vortex}
\end{figure}


\section{TDA and Time Delay Embeddings}
\subsection{Introduction to TDA}
TDA studies the topology of data in an applicable manner. One of the most commonly used methods of TDA is  \textit{persistent homology}. For a brief overview, see Appendix \ref{sec:homology}.

\subsection{Time-delay Embeddings}\label{sec:tdelay}

In this paper, we analyze a one-dimensional time series with a measured variable of low pressure, polar quasizonal winds aggregated into a time series of daily values. 
Applying persistent homology onto this one-dimensional curve would trivially result in no new useful information. 
However, using \textit{Takens' embedding theorem} originally derived for \textit{dynamical systems}, we can construct a higher dimensional point-cloud of our time series that preserves that inherent structure of the original system \cite{takens_detecting_1981}. 
In higher dimensions, we enable more complex topological features to emerge, allowing us to study the \textit{chaotic behavior} or \textit{periodicity} in more depth. 
Takens' embedding, or \textit{time-delay} embedding, can be thought of as sliding a window across a signal. 
Each window is then represented as a point in a space of higher dimension. 

More specifically, take a function $f$ defined on an interval over $\mathbb{R}$. For a chosen $M \in \mathbb{Z}_{>0}$ and $\tau \in \mathbb{R}_{>0}$, the sliding window embedding of $f$ based at $t \in \mathbb{R}$ into $\mathbb{R}^{M+1}$ is the point 
\begin{equation}\label{eqn:SW}
\operatorname{SW}_{M,\tau} f(t) =
\begin{bmatrix}
    f(t) & f(t + \tau) & \ldots & f(t + M\tau)
\end{bmatrix}^T.
\end{equation}
By choosing multiple values for $t$, we obtain a \textit{sliding window point-cloud} for $f$, which will have a non-trivial topology for an appropriate selection of the number of delay coordinates $M$ and sampling frequency $\tau$ \cite{perea_sliding_2015}. 
In the case of weather data, we are able to construct a simplicial complex from this point cloud and ultimately find topological features linked to changes in weather patterns.

\section{Data}
The topological features of certain toy models of different theoretical dynamical systems, such as L63, have been extensively studied in the past. Despite being theoretical, these models reveal extremely topologically persistent features, both from ``looping'' tendencies (like the two butterfly wings of L63) and in their shifts into chaos. Other real-world dynamical systems have been studied using TDA as well, such as the North Atlantic Jet Stream \cite{strommen_topological_2023}.

\subsection{Polar Vortex Data}
The polar vortex is the system of stratospheric winds moving quasi-zonal \cite{schoeberl_structure_1992} around $60$N and $60$S, and can be studied via zonal-mean zonal wind data at low pressures near $60$N or $60$S. Zonal wind is the \textit{u-component} of wind, representing east-west motion. In contrast, the \textit{v-component} represents north-south wind, called meridional wind. Since we are specifically interested in the northern hemisphere, we use data near 60N. We obtained our dataset from the Copernicus Climate Data Store (CDS)\cite{copernicus_climate_change_service_era5_2018}, which is owned and maintained by the European Union.
Upon extraction, the dataset provided a time-series with hourly values for the mean zonal winds. Since we are particularly interested in how the system changes over longer periods of time, this level precision is unnecessary. To preserve computational resources, the time series was averaged into one with a single value for each day. 
We examined data between the beginning of $2015$ and the end of $2020$. We selected this time period because it contained the historically strong and cold polar vortex in the winter of $2015/2016$ \cite{matthias_extraordinarily_2016}, which matched up with Montana temperature data provided by the National Centers for Environmental Information (NCEI) \cite{noauthor_data_nodate}. 

\section{Computation}

\subsection{Persistent Homology}
Persistence diagrams (PD's) are able to uniquely encode topological features computed from persistent homology into a singular diagram. 
We compute persistent homology with the Python package GUDHI \cite{the_gudhi_project_gudhi_2025}.
Our analysis specifically targets the $1$-dimensional homology group $H_1$, which identifies loops and $1$-dimensional holes; see Appendix~\ref{sec:homology} for more details.
We construct our point cloud by leveraging Takens's theorem, highlighted in Section~\ref{sec:tdelay}.
Our chosen delay coordinates are $M = 7$ and $\tau = 1$, creating vectors of $7$ daily observations at a frequency of $\tau = 1$, our fixed window coordinate size being $30$, and sliding window size of $M\tau = 7$ as in Equation \ref{eqn:SW}.
Vietoris--Rips (VR) complexes are then constructed at each stage to compute PD's \cite{edelsbrunner_computational_2010}. 
By shifting the window forward by one observation at each step, we are able to track the temporal changes in the topology of the polar vortex.

\subsection{Visualizing Persistent Homology}
We utilized persistence diagrams and normed persistence landscapes (PL's) (see \cite{bubenik_statistical_2015}) to study the topological features of the polar vortex. 
While PD's are useful in qualitatively assessing how the polar vortex evolves, employing PL's is useful in this case to highlight periods of significant topological activity.
By representing persistence norms as a time series, we can directly compare them to external climatic variables such as temperature or precipitation in regions influenced by the polar vortex. 
Collectively these combined approaches yield a comprehensive view of the topological structure, temporal evolution, and dynamics of the polar vortex.

\begin{figure}[h]
    \centering
    \includegraphics[width=0.62\linewidth]{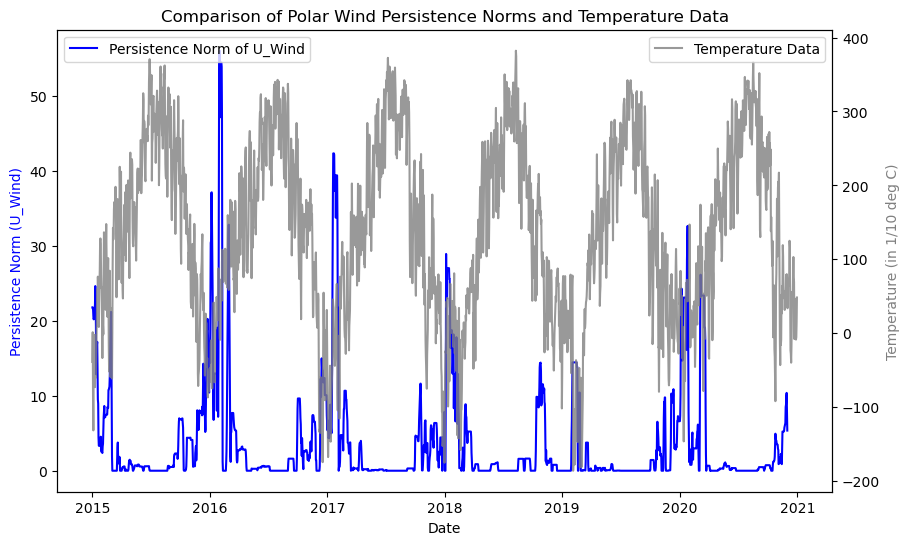}
    \caption{A time series plot of the persistence norms of the polar wind data (blue) compared to aggregate temperature values for Montana, USA, over the period $2015$--$2020$.}
    \label{fig:norms}
\end{figure}

As we may expect, general seasonal changes are very clear. For each year we examined ($2015$--$2020$), the PL's clearly shows the onset and decay of the vortex, since the polar vortex is strongest in winter. This is apparent in Figure~\ref{fig:norms}. 
For each year, as temperatures in Montana start to decrease, we begin to see changes in the topology of the data. 
The significance of these topological features peaks mid-winter and declines as temperatures climb again, usually ending around April. 
The height of the peaks of persistence norms is also noteworthy. Early $2016$ brought the harshest polar vortex event to the United States in, at that time, $86$ years \cite{matthias_extraordinarily_2016}. Our results reflect this, with the highest spike in persistence norms occurring at this time. Each calendar year brings one clear peak, with some more significant than others.

Using PL's is not the only way to visualize these interesting topological features, as we can look directly at the persistence diagrams to gain more qualitative insight as well. Continuing to focus on early $2016$ as an example, Figure~\ref{fig:persistence-a} provides a snapshot of the topological characteristics for a single-point cloud. This PD represents the homology class $H_1$. Since the points are high above the diagonal line, the corresponding topological features differentiate themselves from noise in the data, suggesting that they represent meaningful structures and provide insight into the dynamical system. 

Since this is just a snapshot of our sliding window, we can further understand what is happening by looking at the persistence diagrams in succession. In an animated form, when we looked at how the topological features change, we saw these two features start near the diagonal in the days prior, and then slowly move further from the line until this snapshot on February 1. In the days following, the features move closer together in the diagrams before merging into a single feature, as shown in Figure~\ref{fig:persistence-b}. To view this shift in animated form, see our project repository \cite{majhi_topology_2025}.

\section{Conclusion}
In this paper, we have shown the topological significance of the polar vortex in the northern hemisphere and the influence of the vortex on weather patterns in affected regions. 
In the style of \cite{qin_rapid_2025}, we also entertain the future possibility of using merge tree neural networks for visualizations of this topological activity.

\clearpage


\bibliographystyle{plain}
\bibliography{references}

\appendix

\section{Persistent Homology}\label{sec:homology}
Following the style of \cite{hensel_survey_2021}, a non-empty family of sets $K$ with a collection of non-empty sets $S$ is an \textit{abstract simplicial complex} if both $\{v\} \in S$, for all $v \in K$, and if $\sigma \in S$ and $\tau \subseteq \sigma$, then $\tau \in K$.

The elements of a simplicial complex $K$ are called \textit{simplices}. A $k$-simplex consists of $k+1$ vertices. 
Coefficients are in $\mathbb{Z}_2$, so all elements of $C_p$ are of the form $\sum_j \sigma_j$, for $\sigma_j \in K$. 
Note that $\mathbb{Z}_2$ under addition is often a popular choice since addition also represents \textit{symmetric difference}. 
We use chain groups to express the notion of \textit{boundary}.
Given a simplicial complex $K$, the $p^{\text{th}}$ boundary homomorphism is a function that assigns each simplex $\sigma = \{v_0, \dots, v_p\} \in K$ to its \textit{boundary} $\partial_p \sigma = \sum_{i} \{v_o, \dots, \hat{v}_i, \dots, v_p \}$. 
In this case, $\hat{v}_i$ is the set \textit{not} containing the $i^\text{th}$ vertex, and $\partial_p : C_p \rightarrow C_{p-1}$ is a homomorphism between chain groups. 
For all $p$, $\partial_{p-1} \circ \partial_p = 0$, i.e. boundaries do not have boundaries themselves. This leads to the chain complex $
0 \xrightarrow{\partial_{n+1}} C_n 
\xrightarrow{\partial_n} C_{n-1} 
\xrightarrow{\partial_{n-1}} \dots 
\xrightarrow{\partial_2} C_1 
\xrightarrow{\partial_1} C_0 
\xrightarrow{\partial_0} 0 $.
We define the \textit{cycle group} to be $Z_p = \operatorname{ker}\partial_p$, and the \textit{boundary group} to be $B_p = \operatorname{im}\partial_{p+1}$. 
Note $B_p \subseteq Z_p$ are abelian groups. 
We define the $p^{\text{th}}$ homology group $H_p = Z_p / B_p$, and thus the $p^{\text{th}}$ Betti number $\beta_p = \operatorname{rank}H_p$. 

With this, we can implement \textit{persistent homology}, which computes the appearance (birth) and disappearance (death) of features of a simplicial complex by studying them at different resolutions via a \textit{filtration}. 
A filtration is a nested sequence of subcomplexes (i.e., a simplicial complex $\hat{K_i}$ where $\hat{K_i} \subseteq K$) 
such that $\emptyset = \hat{K_0} \subseteq \hat{K_1} \subseteq \dots \subseteq \hat{K}_{n-1} \subseteq \hat{K_n} = K$. 
Now, since $\hat{K}_i \subseteq \hat{K}_j$ for all $i \leq j$, we get a sequence of homomorphisms connecting the homology groups of each simplicial complex, i.e. $f^{i,j}_{p}: H_p(\hat{K}_i) \rightarrow H_p(\hat{K}_j)$, and thus
$0 = H_p(\hat{K}_0) \xrightarrow{f^{0,1}_{p}} H_p(\hat{K}_1) \xrightarrow{f^{1,2}_{p}} \dots \xrightarrow{f^{n-2,n-1}_{p}} H_p(\hat{K}_{n-1}) \xrightarrow{f^{n-1,n}_{p}} H_p(\hat{K}_n) = H_p(K),$ with $p$ being the dimension of the corresponding homology group. 
Given two indices $i \leq j$, the $p^{\text{th}}$ persistent homology group $H^{i,j}_p$ is defined as
$H^{i,j}_p := Z_p (\hat{K}_i) / \big( B_p (\hat{K}_j) \cap Z_p (\hat{K}_i) \big),$ which contains all the homology classes of $\hat{K}_i$ that are still present in $\hat{K}_j$. 
Thus, we are able to compute the \textit{birth} $(b_i)$ and \textit{death} $(d_i)$ of features across homology groups. 
We then can plot $(b_i, d_i) \in \mathbb R^2$ onto a \textit{persistence diagram} to have a condensed view of the lifetimes of all features across the filtration. 

\section{Additional Figures}
\begin{figure}[hbt]
    \centering
    \begin{subfigure}{.5\textwidth}
        \centering
        \includegraphics[width=0.8\linewidth]{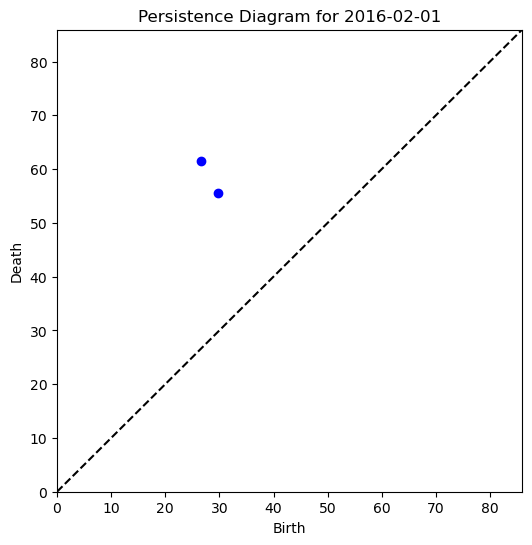}
        \caption{\centering An $H_1$ persistence diagram snapshot from the point cloud of February $1$, $2016$, with two prominent topological features visible}
        \label{fig:persistence-a}
    \end{subfigure}%
    \begin{subfigure}{.5\textwidth}
        \centering
        \includegraphics[width=0.8\linewidth]{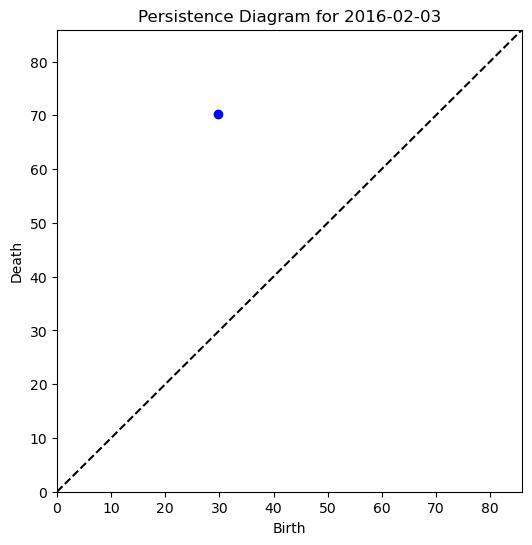}
        \caption{\centering An $H_1$ persistence diagram snapshot from the point cloud of February $03$, $2016$, with now only one topological feature}
        \label{fig:persistence-b}
    \end{subfigure}
    \caption{\label{fig:persistence}}
\end{figure}

\end{document}